\documentclass[a4paper, 11pt]{article}

\usepackage{amsmath,amsfonts,latexsym,theorem%, showkeys
}
\nonstopmode

\numberwithin{equation}{section}

\newtheorem{Theorem}{Theorem}[section]
\newtheorem{Proposition}{Proposition}[section]

\newtheorem{Corollary}{Corollary}[section]
\newenvironment{Proofc}[1]{\smallskip\par\noindent\textsc{#1}\quad}%
  {\hfill$\Box$\bigskip\par}

\newtheorem{Remark}{Remark}[section]
\def\a{\alpha}
\def\b{\beta}
\def\d{\delta}
\def\g{\gamma}

\def\s{\sigma}

\def\t{\tau}
\def\e{\varepsilon}
\newcommand{\Ho}{{\overline H}}
\newcommand{\xo}{{\overline x}}
\newcommand{\po}{{\overline p}}
\newcommand{\Xo}{{\overline X}}
\newcommand{\Co}{{\overline C}}
\newcommand{\Ko}{{\overline K}}
\newcommand{\tr}{\operatorname{\text{tr}}}
 
\newcommand{\R}{{\mathbb R}}
\newcommand{\Rn}{{\mathbb R^n}}
\newcommand{\Rm}{{\mathbb R^m}}
\newcommand{\M}{{\mathbb M}}
\newcommand{\Sym}{{\mathbb S}}

\newcommand{\esssup}{\mathop{\rm ess{\,}sup}}

\begin{document}
\title{Continuous dependence estimates for large time behavior for Bellman-Isaacs equations and applications to the ergodic problem\footnotemark[1] }
\author{Claudio Marchi\footnotemark[2]}
\date{July $28^{\textrm{th}}$, 2010}
%\date{version: \today}
\maketitle

\footnotetext[1]{Work partially supported by
the INDAM-GNAMPA project ``Fenomeni di propagazione di fronti e problemi di
omogeneizzazione''.}
\footnotetext[2]{Dip. di Matematica Pura ed Applicata, Universit\`a di Padova, via Trieste 63, 35121 Padova, Italy ({\tt marchi@math.unipd.it}).}

 \begin{abstract}
This paper concerns continuous dependence estimates for Hamilton-Jacobi-Bellman-Isaacs operators (briefly, HJBI).
For the parabolic Cauchy problem, we establish such an estimate in the whole space $[0,+\infty)\times\Rn$.
Moreover, under some periodicity and ellipticity assumptions, we obtain a similar estimate for the ergodic constant associated to the HJBI operator.
An interesting byproduct of the latter result will be the local uniform convergence for some classes of singular perturbation problems.
\end{abstract}
 
\begin{description}
\item [{\bf MSC 2000}:] 35B25, 35B30, 35J60, 35K55, %35K65,
 49L25, 49N70
.
 \item [{\bf Keywords}:] Continuous dependence estimates, parabolic Hamilton-Jacobi equations, viscosity solutions, ergodic problems, differential games, singular perturbations.

\end{description}

%%%%%%%%%%%%%%%%%%%%%%%%%%%%%%%%
%%        Introduzione        %%
%%%%%%%%%%%%%%%%%%%%%%%%%%%%%%%%

\section{Introduction}\label{intro}

We consider the following Cauchy problem
\begin{equation}\label{P}
\left\{
\begin{array}{ll}
\partial_t u+H(x,Du,D^2u)=0 &\qquad\textrm{in }(0,+\infty)\times\R^n\\
u(0,x)=0&\qquad \textrm{on }\R^n
\end{array}
\right.
\end{equation}
for the Hamilton-Jacobi-Bellman-Isaacs (briefly, HJBI) operator
\begin{equation}\label{minmax}
H(x,p,X)=\min\limits_{\b\in B}\max\limits_{\a\in A}
\left\{
-\tr\left(a(x,\a,\b)X\right)+f(x,\a,\b)\cdot p +\ell(x,\a,\b)
\right\}
\end{equation}
where $\partial_t\equiv \partial/\partial t$,  $Du$ and $D^2u$ stand respectively for the gradient and for the Hessian matrix of the real-valued function~$u=u(t,x)$.

For instance, equations of this kind naturally arise in {\it zero-sum two-persons} stochastic differential games: consider the control system for $s>0$
\begin{equation}\label{control}
dx_s=f(x_s,\a_s,\b_s)+\sqrt 2 \s(x_s,\a_s,\b_s) dW_s,\qquad x_0=x
\end{equation}
where $(\Omega,{\cal F},{\cal P})$ is a probability space endowed with a continuous right filtration $({\cal F}_t)_{0\leq t<+\infty}$ and a $p$-adapted Brownian motion $W_t$.
The control law $\a$ (respectively, $\b$) belongs to the set ${\cal A}$ (resp., ${\cal B}$) of progressively measurable processes which take value in the compact set $A$ (resp., $B$).
The two controls $\a$ and $\b$ are chosen respectively by the first and the second player whose purpose are opposite: the former wants to minimize the {\it cost} functional
\begin{equation}\label{cost}
P(t,x,\a,\b):=\mathbb E_x \int_0^t\ell(x_s,\a_s,\b_s)\, ds
\end{equation}
(here, $\mathbb E_x$ denotes the expectation) while the latter wants to maximize it.
It is well known (see~\cite{FlSo}) that the {\it lower value} function
\begin{equation*}
u(t,x):=\inf\limits_{\a\in \Gamma}\sup\limits_{\b\in {\cal B}}P(t,x,\a[\b],\b)
\end{equation*}
is a viscosity solution to problem~\eqref{P}-\eqref{minmax} with $a={}^T\!\s\s$ where $\Gamma$ stands for the set of admissible {\it strategies} of the first player (namely, nonanticipating maps $\a:{\cal B}\rightarrow {\cal A}$; see~\cite{FlSo}).

\vskip 3mm

This paper is devoted to two main purposes. 
The former purpose is to establish a continuous dependence estimate for problem~\eqref{P}-\eqref{minmax} in the whole space~$[0,+\infty)\times\Rn$; in other words, we want to provide an estimate of $\sup_{\Rn}|u(\cdot,t)-v(\cdot,t)|$ for every $t\in[0,+\infty)$, where $u$ and $v$ are solutions to two problems \eqref{P}-\eqref{minmax} having different coefficients.
The latter purpose of this paper is to establish a continuous dependence estimate for the {\it ergodic constant} associated to HJBI operators~$H$ as in~\eqref{minmax} (see Section~\ref{sect:ergodic} for the precise definition and main properties).
An interesting byproduct of this estimate is the local uniform convergence for some classes of {\it singular perturbation} problems (see Section~\ref{pertsing}).

\vskip 3mm

The continuous dependence estimates for fully nonlinear equations have been widely studied in literature, starting from the paper by Souganidis~\cite{Sou} for first-order equations.
In fact, such estimates play a crucial role in many contexts as error estimates for approximation schemes (see \cite{BJ, DK} and references therein), regularity results (for instance, see \cite{BLM, Bou, JK1}) and rate of convergence for vanishing viscosity methods (see \cite{CGL, Gr, JK1, JK2} and references therein).
In particular, let us recall that Cockburn, Gripenberg and Londen~\cite{CGL} tackled up the continuous dependence estimate for quasi-linear second-order equations with Neumann boundary conditions, while Grinpenberg~\cite{Gr} addressed the case of the Dirichlet boundary data for the same equations.
Afterwards, Jakobsen and Karlsen~\cite{JK1} extended their results to more general classes of equations (see also \cite{JK2} for elliptic problems).
Furthermore, Jakobsen and Georgelin~\cite{JG} extended the previous results to problems with more general boundary conditions and domains.

The first main purpose of this paper is to establish a continuous dependence estimate for problem~\eqref{P}-\eqref{minmax} with the following three features:
the estimate holds in the whole space $[0,+\infty)\times\Rn$,  the dependence on the $\mathbb L^\infty$-distance between the coefficients is explicit, the constants can be explicitly characterized.

As one may expect, it turns out that our estimate increases linearly with $t$.
A similar estimate could be obtained by an easy application of the Comparison Principle provided that some bound on the $C^2$-norm of the solution is available.
Unfortunately, this is not the case for operators in form~\eqref{minmax}.
In fact, our approach is based on the Comparison Principle techniques for viscosity solutions (see~\cite{CIL}): doubling the variables and adding a penalization term.
Let us observe that this approach does not require the non-degeneracy of the operator~$H$; actually, we shall also apply our result to some degenerate problems.

\vskip 3mm

In the {\it ergodic problem} for the operator $H$, we seek a pair $(v,U)$, with $v\in C(\Rn)$ and $U\in \R$ (that is, $v$ is a real-valued  function while $U$ is a constant) which, in viscosity sense, satisfy
\begin{equation}\label{E1}
H(x,Dv,D^2v)=U\qquad\textrm{in }\R^n.
\end{equation}
This problem has been widely studied in literature, especially in connection with homogenization or singular perturbation problems (see~\cite{AB2, AB8, BDLS, CSW, Ev2, LS2} and references therein), with long-time behavior of solutions to parabolic equations (for instance, see \cite{AL, BaDa, BaSo}) and with dynamical systems in a torus (see~\cite{AA, CFS, Si}).

It is well known (see \cite{ AB8,AL}) that, under some periodicity and non-degeneracy  assumptions, there exists exactly one value $U\in\R$ (called {\it ergodic constant}) such that equation~\eqref{E1} admits at least one bounded solution (which is periodic and unique up to a constant).
In \cite{Ev2}, Evans obtained a continuous dependence estimate for the ergodic constant for operators which are Lipschitz continuous in the variable $x$ uniformly with respect to $(p,X)$.
Afterwards, Alvarez and Bardi~\cite{AB8}, extended his result to operators $H$ as in~\eqref{minmax} provided that the dispersion matrix is left unchanged (see also \cite[Section 6]{AB8} for possibly degenerate equations).

The second main purpose of this paper is to obtain a continuous dependence estimate for the ergodic constant of HJBI operators~\eqref{minmax} only under the periodicity and the non-degeneracy assumption namely, we shall also consider equations with different dispersion matrices.
An interesting byproduct will be the local uniform convergence for some classes of singular perturbation problems for HJBI operators.

\vskip 3mm

In conclusion, the aim of this paper is threefold: a continuous dependence estimate for problem~\eqref{P}-\eqref{minmax} in the whole space $[0,+\infty)\times\Rn$, a similar estimate for the ergodic constant for HJBI operators \eqref{minmax} and a local uniform convergence for singular perturbation problems.

This paper is organized as follows: in the rest of this section, we provide some notations and list the standing assumptions.
Section~\ref{results} contains the continuous dependence estimate for the parabolic Cauchy problem~\eqref{P} and its application to some degenerate problems as well.
Section~\ref{sect:ergodic} concerns the continuous dependence estimate for the ergodic constant in~\eqref{E1}; section~\ref{pertsing} is devoted to illustrate how to derive the local uniform convergence for singular perturbation problems.

%%%$%%%%%%%%%%%%%%%%%%%%%%%%%%%%%%%%%%
%% notazioni e standing assumptions %%
%%%%%%%%%%%%%%%%%%%%%%%%%%%%%%%%%%%%%%

\subsection{Notations and standing assumptions}\label{assumption}
\textbf{Notations: } We define $\mathbb M^{n,p}$ and $\Sym^n$ respectively as the set of
$n\times p$ real matrices and the set of $n\times n$ symmetric matrices. The latter is endowed with the Euclidean norm and with the usual order, namely: for $X=(X_{ij})_{i,j=1,\dots,n}\in \Sym^n$, $\|X\|:=(\sum_{i,j=1}^n X^2_{ij})^{1/2}$ and, for $X,Y\in\Sym^n$, we shall write $X\geq Y$, if $X-Y$ is a semi-definite positive matrix.

For every positive $t$, we set $Q_t:=[0,t)\times\Rn$ and $Q_\infty:=[0,+\infty)\times\Rn$.

For every real-valued function $h$, we set $\|h\|_\infty:=\esssup|h(y)|$; for $\g\in(0,1]$, we use the $\g$-H\"older norm: $|h|_\g:=\sup_{y\neq x} \frac{|h(y)-h(x)|}{|y-x|^\g}$.
Moreover, $J^{2,+}h(\xi)$ and $J^{2,-}h(\xi)$ stand respectively for the second-order superjet and subjet of $h$ at the point $\xi$ (see \cite{CIL} for the precise definition and main properties).
A real function $\omega$ is said a {\it modulus of continuity} whenever it is a nonnegative continuous non-decreasing real function on $[0,+\infty)$ with $\omega(0)= 0$.

\vskip 1mm

\noindent\textbf{Standing assumptions:}
For the operator~$H$ in~\eqref{minmax}, the following assumptions will hold throughout this paper
\begin{itemize}
\item[($A1$)] $A$ and $B$ are two compact metric spaces.
\item[($A2$)] $a={}^T\!\s\s$. The functions $\s$, $f$ and~$\ell$ are bounded continuous functions in $\Rn \times A \times B$ with value respectively in $\M^{n,p}$, $\Rn$, and $\R$; namely, for some $C>0$, there holds: $\|\s\|_\infty,\, \|f\|_\infty,\, \|\ell\|_\infty\leq C$.
\item[($A3$)] The drift vectors $f$ and the dispersion matrix $\s$ are Lipschitz continuous in $x$ uniformly in $(\a,\b)$, namely: for some positive constant $C_\phi$, every $\phi=\s, f$ satisfies
\begin{equation*}
\left|\phi(x,\a,\b)-\phi(y,\a,\b)\right|\leq C_\phi |x-y|\qquad \forall x,y\in\Rn,\, \forall (\a,\b)\in A\times B.
\end{equation*}
The running cost~$\ell$ is uniformly continuous in $x$ uniformly in $(\a,\b)$, namely: there exists a modulus of continuity $\omega$ such that
\begin{equation*}
\left|\ell(x,\a,\b)-\ell(y,\a,\b)\right|\leq \omega(|x-y|)\qquad \forall x,y\in\Rn,\, \forall (\a,\b)\in A\times B.
\end{equation*}
\end{itemize}

%%%%%%%%%%%%%%%%%%%%%%%%%%%%%%%%
%% cde delproblema parabolico %%
%%%%%%%%%%%%%%%%%%%%%%%%%%%%%%%%

\section{Estimate for the parabolic Cauchy problem}\label{results}

For $i=1,2$, consider the parabolic Cauchy problems 
\begin{equation}\label{Pi}\tag{Pi}
\left\{
\begin{array}{ll}
\partial_t u_i+\min\limits_{\b\in B}\max\limits_{\a\in A}\left\{
-\tr\left(a_i(x,\a,\b)D^2u_i\right)+f_i(x,\a,\b)\cdot Du_i +\ell_i(x,\a,\b)
\right\}=0& \textrm{ in }Q_\infty\\
u_i(0,x)=0& \textrm{ on }\R^n.
\end{array}
\right.
\end{equation}
where the coefficients fulfill our standing assumptions ($A1$)-($A3$).
The main purpose of this section is to provide an estimate of $\|u_1(t,\cdot)-u_2(t,\cdot)\|_\infty$ for every $t\in[0,+\infty)$. In section~\ref{Example} we shall apply this estimate to some degenerate problems.

\begin{Remark}\label{exi-ui}
By standard viscosity theory (for instance, see \cite{CIL}), assumptions ($A1$)-($A3$) guarantee that the Comparison Principle applies to problem~\eqref{Pi}; whence, by the Perron method, one can easily deduce that \eqref{Pi} admits exactly one solution $u_i\in C(Q_\infty)$ with
\begin{equation}\label{exi-uii}
|u_i(t,x)|\leq tC,\qquad \forall(t,x)\in Q_\infty
\end{equation}
where $C$ is the constant introduced in assumption~($A2$).
\end{Remark}

%
%Teorema cde problema parabolico
%

\begin{Theorem}\label{thm:dipcnt}
Let $u_i$ be the unique solution to problem~\eqref{Pi} which satisfies the bound~\eqref{exi-uii} ($i=1,2$).
Furthermore, let us assume that, for some $\g\in(0,1]$, $u_i(t,\cdot)$ is $\g$-H\"older continuous uniformly in $t$, namely: for some $C_H>0$, there holds 
\begin{equation}\label{A100}
|u_i(t,\cdot)|_\g\leq C_H,\qquad \forall t\in[0,+\infty),\,i=1,2.
\end{equation}
Then, there exist a positive constant $M$ such that, for every $(x,t)\in Q_\infty$, there holds
\begin{multline*}
\left|u_1(t,x)-u_2(t,x)\right|\leq
t M \left[\max\limits_{x,\a,\b}\left\|\s_1-\s_2\right\|^\g
+\max\limits_{x,\a,\b}|f_1-f_2|^{\g/2}+\max\limits_{x,\a,\b}|\ell_1-\ell_2|\right.\\ \left.
+\omega\left(C(\max\limits_{x,\a,\b}\left\|\s_1-\s_2\right\|
+\max\limits_{x,\a,\b}|f_1-f_2|^{1/2})\right)\right].
\end{multline*}
\end{Theorem}

%
%Dimostrazione teo cde parabolico
%

\begin{Proofc}{Proof of Theorem~\ref{thm:dipcnt}}
We shall argue using some techniques introduced in~\cite{CGL,JK1}.
We fix $t>0$ and, for every $\eta\in (1,+\infty)$ and $\e\in(0,1)$, we introduce $E_t:=[0,t)\times\Rn\times\Rn$ and
\begin{equation}\label{sigma}
s_t:=\sup\limits_{E_t}\left\{u^1(\t,x)-u^2(\t,y)
-\left(\frac \eta 2 |x-y|^2+\frac \e2(|x|^2+|y|^2)+\frac {\e}{t-\tau}\right)\right\}.
\end{equation}
Our purpose is to establish an upper bound for $s_t$. To this end, without any loss of generality, we can assume $s_t>0$.
For $\d\in(0,1)$, define
\begin{equation}\label{psi}
\psi(\t,x,y):=u^1(\t,x)-u^2(\t,y)-\frac{\d s_t\t}{t}-
\left(\frac \eta 2 |x-y|^2+\frac \e2(|x|^2+|y|^2)
+\frac {\e}{t-\tau}\right).
\end{equation}
Let us observe that definition~\eqref{sigma} entails
\begin{eqnarray}
\sup\limits_{E_t}\psi 
&\geq& \sup\limits_{E_t}\left\{u^1(\t,x)-u^2(\t,y)-
\left(\frac \eta 2 |x-y|^2+\frac \e2 (|x|^2+|y|^2)+\frac {\e}{t-\tau}\right)
\right\}\notag\\
&&\hskip 10mm
-\sup\limits_{\t\in[0,t)}\left\{\frac{\d s_t\t}{t}\right\}\notag\\ \label{maxpos}
&\geq& (1-\d) s_t>0.
\end{eqnarray}
Since the functions $u^1$ and $u^2$ are bounded in~$Q_t$ and since the function $\psi$ tends to $-\infty$ as $\t\to t^-$, we deduce that there exists a point~$(\t_0,x_0,y_0)\in E_t$ where the function~$\psi$ attains its global maximum, namely
\begin{equation}\label{psi0>0}
\psi(\t_0,x_0,y_0)=\sup\limits_{E_t}\psi\geq0
\end{equation}
where last inequality is due to relation~\eqref{maxpos}.

Let us now claim that, for  $\Co:=(2C_H)^{1/(2-\g)}$, there holds
\begin{equation}\label{stima:xy}
|x_0-y_0|\leq \Co \eta^{-1/(2-\g)},\qquad \e\left(|x_0|^2+|y_0|^2\right)\leq 4Ct
\end{equation}
where $C_H$, $\g$ and $C$ are the constants introduced respectively in assumption~\eqref{A100} and ($A2$); in particular, let us emphasize that $\Co$ is independent of $t$.
Actually, in order to prove the former estimate, we observe that inequality $\psi(\t_0,x_0,x_0)+\psi(\t_0,y_0,y_0)\leq 2\psi(\t_0,x_0,y_0)$ and assumption~\eqref{A100} give
\begin{equation*}
\begin{array}{rcl}
\eta|x_0-y_0|^2&\leq& [u^1(\t_0,x_0)-u^1(\t_0,y_0)]+[u^2(\t_0,x_0)-u^2(\t_0,y_0)]\\
&\leq& 2C_H|x_0-y_0|^\g.
\end{array}
\end{equation*}
Let us now prove the latter estimate in~\eqref{stima:xy}: by estimates~\eqref{exi-uii} and~\eqref{psi0>0}, we infer
\begin{equation*}
\e\left(|x_0|^2+|y_0|^2\right)\leq 2 u_1(\t_0,x_0)-2 u_2(\t_0,y_0)\leq 4Ct.
\end{equation*}
Hence, the proof of estimates~\eqref{stima:xy} is accomplished.

We introduce the test function
\begin{equation}\label{phi}
\phi(\t,x,y):=
\frac{\d s_t\t}{t}+\frac \eta 2 |x-y|^2+\frac\e2(|x|^2+|y|^2)+\frac {\e}{t-\tau}
\end{equation}
and we invoke \cite[Theorem 8.3]{CIL}: for every $\nu>0$, there exist values $a,b\in\R$ and matrices $X,Y\in\Sym^n$ such that
\begin{eqnarray}
&\label{cil1}
\left(a,D_x\phi(\t_0,x_0,y_0),X\right)\in J^{2,+}u^1(\t_0,x_0),\quad
\left(b,D_y\phi(\t_0,x_0,y_0),Y\right)\in J^{2,-}u^2(\t_0,y_0),\\
&a-b=\partial_\tau\phi(\t_0,x_0,y_0)\equiv \frac{\d s_t}t+\frac {\e}{(t-\tau_0)^2}&\label{cil2} \\
&\left(\begin{array}{cc} X&0\\0&-Y\end{array}\right)\leq
\Theta+\nu \Theta^2, & \label{cil3}
\end{eqnarray}
where $\Theta:=\eta\left(\begin{array}{cc} I&-I\\-I&I\end{array}\right) + \e\left(\begin{array}{cc} I&0\\0&I\end{array}\right)$.
%
%Trick di Ishii '89
%
From the last inequality, one can deduce that, for every $(\a,\b)\in A\times B$, there holds
\begin{multline}\label{Ishii}
\tr\left(a_1(x_0,\a,\b)X\right)-
\tr\left(a_2(y_0,\a,\b)Y\right)\leq \eta\left\|\s_1(x_0,\a,\b)-\s_2(y_0,\a,\b)\right\|^2 \\
+2\e C^2+\nu \tr\left(\Sigma\Theta^2\right)
\end{multline}
with
$$
\Sigma:=\left(\begin{array}{lc}
{}^T\!\s_1(x_0,\a,\b)\s_1(x_0,\a,\b)&{}^T\!\s_1(x_0,\a,\b)\s_2(y_0,\a,\b)\\
{}^T\!\s_2(y_0,\a,\b)\s_1(x_0,\a,\b)&{}^T\!\s_2(y_0,\a,\b)\s_2(y_0,\a,\b)
\end{array}\right).
$$
In order to prove this inequality, we shall use the arguments by Ishii~\cite{Is89}.
Multiplying relation \eqref{cil3} by matrix~$\Sigma$ (which is symmetric and nonnegative definite)
and evaluating the trace, we obtain
\begin{align*}
&\tr\left({}^T\!\s_1(x_0,\a,\b)\s_1(x_0,\a,\b)X-
{}^T\!\s_2(y_0,\a,\b)\s_2(y_0,\a,\b)Y\right)\\
&\hskip 10mm \leq 
\eta\tr\left[{}^T\!\left(\s_1(x_0,\a,\b)-\s_2(y_0,\a,\b)\right)\left(\s_1(x_0,\a,\b)-\s_2(y_0,\a,\b)\right)\right]\\
&\hskip20mm
+\e \tr\left({}^T\!\s_1(x_0,\a,\b)\s_1(x_0,\a,\b)\right)
+\e \tr\left({}^T\!\s_2(y_0,\a,\b)\s_2(y_0,\a,\b)\right)
+\nu \tr\left(\Sigma\Theta^2\right)
;
\end{align*}
therefore, by assumption~($A2$), relation~\eqref{Ishii} easily follows.
%fine dim trick di Ishii

Since $u^1$ is a subsolution to problem~\eqref{Pi}  with $i=1$, the former relation in~\eqref{cil1} yields
\begin{eqnarray*}
0&\geq& a+\min\limits_{\b\in B}\max\limits_{\a\in A}\left\{
-\tr\left(a_1(x_0,\a,\b)X\right)+f_1(x_0,\a,\b)\cdot D_x\phi +\ell_1(x_0,\a,\b)
\right\}\\
&\geq& b+\min\limits_{\b\in B}\max\limits_{\a\in A}\left\{
-\tr\left(a_2(y_0,\a,\b)Y\right)+ f_1(x_0,\a,\b)\cdot (\eta(x_0-y_0)+\e x_0) +\ell_1(x_0,\a,\b)\right\}\\
&&\hskip3mm -\eta\max\limits_{\a,\b}\left\|\s_1(x_0,\a,\b)-\s_2(y_0,\a,\b)\right\|^2 -2\e C^2 -\nu\max\limits_{\a,\b} \tr\left(\Sigma \Theta^2\right)
+\frac{\d s_t}t+\frac {\e}{(t-\tau_0)^2}
\end{eqnarray*}
where the last inequality is due to the definition of $\phi$ \eqref{phi} and to relations~\eqref{cil2} and~\eqref{Ishii}.
Since $u^2$ is  a supersolution to equation~\eqref{Pi}  with $i=2$, by assumption ($A2$), last inequality entails the following upper bound for $s_t$:
\begin{eqnarray*}
\frac{\d s_t}t+\frac {\e}{(t-\tau_0)^2}&\leq&\eta\max\limits_{\a,\b}\left\|\s_1(x_0,\a,\b)-\s_2(y_0,\a,\b)\right\|^2
+2\e C^2 +\nu \max\limits_{\a,\b}\tr\left(\Sigma \Theta^2\right)\\
&&\hskip 5mm+\e C \left(|x_0|+ |y_0|\right)
+\eta|x_0-y_0|\max\limits_{\a,\b}\left|f_1(x_0,\a,\b)-f_2(y_0,\a,\b)\right|\\
&&\hskip 5mm+
\max\limits_{\a,\b}\left|\ell_1(x_0,\a,\b)-\ell_2(y_0,\a,\b)\right|.
\end{eqnarray*}
Owing to the definition of $s_t$ in~\eqref{sigma}, we deduce
\begin{equation*}
u^1(\t,x)-u^2(\t,y)-\frac \eta 2 |x-y|^2-\frac \e2\left(|x|^2+|y|^2\right)\leq
 s_t+\frac {\e}{t-\tau}\qquad\forall (\t,x,y)\in E_t.
\end{equation*}
Taking into account the last two inequalities, for every $(\t,x)\in Q_t$, we infer
\begin{align*}
&u^1(\t,x)-u^2(\t,x)\\
&\hskip10mm \leq
\frac t\d\left[
\eta\max\limits_{\a,\b}\left\|\s_1(x_0,\a,\b)-\s_2(y_0,\a,\b)\right\|^2
+2\e C^2 +\nu \max\limits_{\a,\b}\tr\left(\Sigma\Theta^2\right)
\right.\\
&\hskip20mm +\eta|x_0-y_0|\max\limits_{\a,\b}\left|f_1(x_0,\a,\b)-f_2(y_0,\a,\b)\right|
+\e C \left(|x_0|+ |y_0|\right)
\\
&\hskip20mm\left.+\max\limits_{\a,\b}\left|\ell_1(x_0,\a,\b)-\ell_2(y_0,\a,\b)\right|
\right] +\frac {\e}{t-\tau}+ \e |x|^2.
\end{align*}
By the regularity of the coefficients (see assumption~ ($A3$)) and estimate~\eqref{stima:xy}, for $\tilde C:=2C_\s^2\Co^2+2+C_f\Co^2+\Co$, we have
\begin{align*}
&u^1(\t,x)-u^2(\t,x)\\
&\hskip 3mm \leq
\frac t\d\left[
 2\eta\left(C^2_\s|x_0-y_0|^2+\max\limits_{x,\a,\b}\left\|\s_1-\s_2\right\|^2
\right)
+\eta C_f|x_0-y_0|^2+\eta|x_0-y_0|\max\limits_{x,\a,\b}|f_1-f_2|
\right.\\&\hskip8mm\left.
+\omega(|x_0-y_0|)+\max\limits_{x,\a,\b}|\ell_1-\ell_2|\right]
+\e \left[\frac 1{t-\tau} +\frac{Ct}{\d} (|x_0|+|y_0|+2C)+|x|^2\right]
\\&\hskip8mm
+\nu \max\limits_{\a,\b}\tr\left(\Sigma\Theta^2\right)
\\&\hskip3mm\leq
\frac t\d\tilde C\left[\eta^{-\g/(2-\g)} +\eta\max\limits_{x,\a,\b}\left\|\s_1-\s_2\right\|^2 +\eta^{(1-\g)/(2-\g)}\max\limits_{x,\a,\b}|f_1-f_2|\right]
\\&\hskip8mm
+\frac t\d\left[\max\limits_{x,\a,\b}|\ell_1-\ell_2|+\omega(\Co\eta^{-1/(2-\g)})\right]
+\e \left[\frac 1{t-\tau} +\frac{Ct}{\d} (|x_0|+|y_0|+2C)+|x|^2\right]\\&\hskip8mm
+\nu \max\limits_{\a,\b}\tr\left(\Sigma\Theta^2\right).
\end{align*}
Letting $\nu\to 0^+$ and afterwards $\e\to 0^+$, by estimate~\ref{stima:xy}, we infer
\begin{multline*}
u^1(\t,x)-u^2(\t,x)\leq
\frac t\d\tilde C\left[\eta^{-\g/(2-\g)} +\eta\max\limits_{x,\a,\b}\left\|\s_1-\s_2\right\|^2
+\eta^{(1-\g)/(2-\g)}\max\limits_{x,\a,\b}|f_1-f_2|
\right]\\+\frac t\d \left[\max\limits_{x,\a,\b}|\ell_1-\ell_2|+
\omega(\Co\eta^{-1/(2-\g)})\right]
\end{multline*}
Letting $\d\to 1^-$ and afterwards $\t\to t^-$, by the continuity of the functions $u^1$ and $u^2$, for every $x\in\Rn$, we deduce
\begin{multline*}
u^1(t,x)-u^2(t,x)\leq
 t\tilde C\left[\eta^{-\g/(2-\g)} +\eta\max\limits_{x,\a,\b}\left\|\s_1-\s_2\right\|^2
+\eta^{(1-\g)/(2-\g)}\max\limits_{x,\a,\b}|f_1-f_2|
\right]\\+ t\left[\max\limits_{x,\a,\b}|\ell_1-\ell_2|
+\omega(\Co\eta^{-1/(2-\g)})\right]
\end{multline*}
Since $\eta$ belongs to $[1,+\infty)$, we infer
\begin{multline*}
u^1(t,x)-u^2(t,x)\leq
t\tilde C \left[\eta^{-\g/(2-\g)} +\eta\left(\max\limits_{x,\a,\b}\left\|\s_1-\s_2\right\|^2+\max\limits_{x,\a,\b}|f_1-f_2|\right)\right]
\\+t\left[\max\limits_{x,\a,\b}|\ell_1-\ell_2|+\omega(\Co\eta^{-1/(2-\g)})
\right]
\end{multline*}
Observe that, for $r\in(0,1)$, the minimum of $h(s):=r s +s^{-\g/(2-\g)}$ in $[1,+\infty)$ is less or equal to $2r^{\g/2}$ (this value is attained in $s=r^{-(2-\g)/2}$).
Therefore, choosing $\eta=[\max\limits_{x,\a,\b}\left\|\s_1-\s_2\right\|^2+\max\limits_{x,\a,\b}|f_1-f_2|]^{-(2-\g)/2}$, we conclude
\begin{eqnarray*}
u^1(t,x)-u^2(t,x)&\leq&
2t\tilde C \left(\max\limits_{x,\a,\b}\left\|\s_1-\s_2\right\|^2
+\max\limits_{x,\a,\b}|f_1-f_2|\right)^{\g/2}\\ &&\hskip 10mm
+t\left[\max\limits_{x,\a,\b}|\ell_1-\ell_2|+\omega\left(\Co(\max\limits_{x,\a,\b}\left\|\s_1-\s_2\right\|^2+\max\limits_{x,\a,\b}|f_1-f_2|)^{1/2}\right)\right]\\ &\leq&
t\tilde C \left[\max\limits_{x,\a,\b}\left\|\s_1-\s_2\right\|^\g
+\max\limits_{x,\a,\b}|f_1-f_2|^{\g/2}\right]\\
&&\hskip 10mm
+t\left[\max\limits_{x,\a,\b}|\ell_1-\ell_2|+\omega\left(\Co(\max\limits_{x,\a,\b}\left\|\s_1-\s_2\right\|+\max\limits_{x,\a,\b}|f_1-f_2|^{1/2})\right)\right]
\end{eqnarray*}
for every $x\in\Rn$.
Owing to the arbitrariness of the value $t$, one side of the inequality in our statement is completely proved. Being similar, the proof of the other one is omitted.
\end{Proofc}

%%%%%%%%%%%%%%%%%%%%%%%%%%%%%%%%
%%   EX: problema degenere    %%
%%%%%%%%%%%%%%%%%%%%%%%%%%%%%%%%

\subsection{Example: a degenerate parabolic problem}\label{Example}
This Section is devoted to illustrate an application of Theorem~\ref{thm:dipcnt} to some classes of parabolic Cauchy problems for degenerate HJBI operators.
\begin{Corollary}\label{Prp:coer}
Assume that, besides our standing assumptions, for some $\nu>0$, there holds
\begin{equation}\label{coer}
\min\limits_{\b\in B}\max\limits_{\a\in A}
\left\{
-\tr\left(a_i(x,\a,\b)X\right)+f_i(x,\a,\b)\cdot p +\ell_i(x,\a,\b)
\right\}
\geq \nu |p|-C
\end{equation}
for every $(x,p,X)\in\Rn\times\Rn\times\Sym^n$ ($i=1,2$).
Then, there exists $M>0$ such that, for every $(t,x)\in Q_\infty$, there holds
\begin{multline*}
\left|u_1(t,x)-u_2(t,x)\right|\leq
t M \left[\max\limits_{x,\a,\b}\left\|\s_1-\s_2\right\|
+\max\limits_{x,\a,\b}|f_1-f_2|^{1/2}+\max\limits_{x,\a,\b}|\ell_1-\ell_2|\right.\\ \left.
+\omega\left(C(\max\limits_{x,\a,\b}\left\|\s_1-\s_2\right\|
+\max\limits_{x,\a,\b}|f_1-f_2|^{1/2})\right)\right],
\end{multline*}
where $u_1$ and $u_2$ are respectively the solution to~\eqref{Pi} with $i=1$ and $i=2$.
\end{Corollary}
\begin{Remark}
Relation~\eqref{coer} is fulfilled provided that there exists $A_i'\subset A$ such that
\begin{equation*}
\sigma_i(x,\a,\b)=0\qquad\forall \a\in A_i',\qquad
B(0,\nu)\subset\overline{\textrm{conv}}\{f_i(x,\a,\b)\mid \a\in A_i'\}
\end{equation*}
for every $x\in\Rn$, $\b\in B$ (here, $B(0,\nu)$ stands for the ball centered in $0$ with radius~$\nu$ while ${\textrm{conv}}{\cal A}$ is the convex hull of ${\cal A}\subset \Rn$).
\end{Remark}
\begin{Proofc}{Proof of Corollary~\ref{Prp:coer}}
A straightforward application of Theorem~\ref{thm:dipcnt} yields the statement provided that the functions $u_1$ and $u_2$ satisfy condition~\eqref{A100} with $\g=1$.
Let us prove this property by using some arguments of~\cite[Theorem~II.1]{AL}.
Assume that there holds
\begin{equation}\label{lip-t}
|u_i(t,x)-u_i(t+h,x)|\leq C h, \qquad \forall (t,x)\in Q_\infty,\, h>0,\,i=1,2.
\end{equation}
In this case, relations~\eqref{coer} and~\eqref{lip-t} guarantee in viscosity sense
\begin{equation*}
C\geq H(x,Du_i(\cdot,t),D^2u_i(\cdot,t))\geq \nu |Du_i(\cdot,t)|-C,  \qquad\textrm{in }\Rn
\end{equation*}
for all $t\in[0,+\infty)$, $i=1,2$.
In particular, we have: $|Du_i|\leq 2C\nu^{-1}$, which amounts to \eqref{A100} with $\g=1$.

In conclusion, let us prove inequality~\eqref{lip-t}.
By estimate~\eqref{exi-uii}, we infer that the functions~$u_i(t+h,x)\pm Ch$ are respectively a super and a subsolution to~\eqref{Pi}.
Applying the Comparison Principle, we accomplish the proof of estimate~\eqref{lip-t}.
\end{Proofc}

%%%%%%%%%%%%%%%%%%%%%%%%%%%%%%%%
%%  cde delproblema ergodico  %%
%%%%%%%%%%%%%%%%%%%%%%%%%%%%%%%%

\section{Estimate for the ergodic problem}\label{sect:ergodic}

This section is devoted to provide a continuous dependence estimate for the ergodic constant associated to the HJBI operator \eqref{minmax}.
Let us recall that, in the ergodic problem, we seek a constant~$U$ such that the equation
\begin{equation}\label{E}
H(x,Dv,D^2v)=U\qquad \textrm{in }\Rn
\end{equation}
admits at least one solution~$v$.
For $\d>0$, let us also introduce the {\it approximated} equation
\begin{equation}\label{Eid}
\d w_\d+H\left(x,Dw_\d,D^2w_\d\right)=0\qquad \textrm{in }\R^n.
\end{equation}

Beside our standing assumptions, throughout this section, the operator~$H$ also fulfills
\begin{itemize}
\item[($A4$)] Periodicity: the functions $\s$, $f$ and~$\ell$ are $\mathbb Z^n$-periodic in $x$;
\item[($A5$)] Non-degeneracy: there exists a positive constant $\nu$ such that:
\begin{equation*}
a(x,\a,\b)\geq \nu I, \qquad\forall(x,\a,\b)\in\Rn\times A\times B.
\end{equation*}
\end{itemize}

For later use, in the following Proposition, we shall collect several known properties of the ergodic problem.
\begin{Proposition}\label{Prp:erg}
Under assumptions ($A1$)-($A5$), the following properties hold:
\begin{itemize}
\item[$i$)] There exists exactly one constant $U$ such that equation~\eqref{E} admits a bounded continuous (and periodic) solution~$v$. Moreover, $v$ is unique up to an additive constant.
\item[$ii$)] Let $u$ be the solution to the Cauchy problem~\eqref{P}; then, as $t\to+\infty$, $u(t,x)/t$ converges to the ergodic constant~$U$ of equation~\eqref{E} uniformly in $x$.
\item[$iii$)] The approximated equation~\eqref{Eid} admits exactly one bounded continuous solution~$w_\d$: $\d \|w_\d\|_\infty\leq \max_{x,\a,\b}|\ell|$. Moreover, as $\d\to 0^+$, $\d w_\d$ and $(w_\d-w_\d(0))$ converge respectively to the ergodic constant~$U$ and to the solution~$v$ of~\eqref{E} with $v(0)=0$.
\item[$iv$)] There exist two constants $\kappa\in(0,1]$ and $K>0$, both depending only on the parameters of our assumptions (that is, independent of $\d$) such that there holds
\begin{equation*}
\|w_\d -w_\d(0)\|_{C^{1,\kappa}}\leq K\left(1+\max\limits_{x,\a,\b}|\ell|\right).
\end{equation*}
\end{itemize}
\end{Proposition}
\begin{Proofc}{Proof of Proposition~\ref{Prp:erg}}

The proof of points~($i$), ($ii$) and~($iii$) can be found in~\cite[Theorem~4.1]{AB8} (see also~\cite[Theorem II.2]{AL} for operators of Bellman type). In fact, the statement of points~($ii$) and~($iii$) are equivalent (see \cite[Theorem 4]{AB2}, \cite[Proposition~2.2]{AB8} and also~\cite[Proposition~VI.1]{AL} for Bellman operators) while the first part of point~($i$) is only a sufficient condition for them (see~\cite[Proposition~7.2]{AB8}).

The proof of point~($iv$) can be easily obtained adapting to HJBI equations the arguments introduced by Arisawa and Lions~\cite[Theorem~II.2]{AL} (see also \cite[Theorem~4.1]{AB8} for a similar result). Finally, we refer the reader to \cite[proof of Theorem~4.1]{AB8} for the special form of the right-hand side.
\end{Proofc}

For $i=1,2$, consider the ergodic problems
\begin{equation}\label{Ei}\tag{Ei}
\min\limits_{\b\in B}\max\limits_{\a\in A}\left\{
-\tr\left(a_i(x,\a,\b)D^2v^i\right)+f_i(x,\a,\b)\cdot Dv^i +\ell_i(x,\a,\b)
\right\}=U^i\quad \textrm{in }\R^n
\end{equation}
where the coefficients fulfill assumptions ($A1$)-($A5$).

\begin{Theorem}\label{thm:erg}
Let $U^i$ be the unique ergodic constant for problem~\eqref{Ei} ($i=1,2$).
Then, there exist a positive constant $\tilde M$ such that there holds
\begin{equation*}
\left|U^1-U^2\right|\leq
\tilde M \left(\max\limits_{x,\a,\b}\left\|\s_1-\s_2\right\|
+\max\limits_{x,\a,\b}|f_1-f_2|\right)
 +\omega(\max\limits_{x,\a,\b}\left\|\s_1-\s_2\right\|)
+\max\limits_{x,\a,\b}|\ell_1-\ell_2|.
\end{equation*}
\end{Theorem}

\begin{Proofc}{Proof of Theorem~\ref{thm:erg}}
Let $u_i$ be the solution to problem~\eqref{Pi} ($i=1,2$).
By Proposition~\ref{Prp:erg}-($ii$), the statement follows from a straightforward application of Theorem~\ref{thm:dipcnt} provided  that the solutions $u_1$ and $u_2$ fulfill condition~\eqref{A100} with $\g=1$.
In order to prove this fact, we denote by $v^i$ the unique bounded solution to equation~\eqref{Ei} with $v^i(0)=0$ and we introduce the function $w^i(t,x):=u_i(t,x)+U^it$, which is the unique solution to the Cauchy problem
\begin{equation}\label{ei}
\left\{
\begin{array}{ll}
\partial_t w^i+\min\limits_{\b\in B}\max\limits_{\a\in A}\left\{
-\tr\left(a_i D^2w^i\right)+f_i\cdot Dw^i +\ell_i\right\}= U^i
&\qquad \textrm{in }Q_\infty\\
w^i(0,x)=0&\qquad \textrm{on }\R^n.
\end{array}
\right.
\end{equation}
Let us prove that $w^i$ is bounded in $Q_\infty$ arguing as in~\cite[Theorem II.1]{AL}.
For $k:=\|v^i\|_\infty$, the functions $v^i(x)-U^it \pm k$ are respectively a super- and a subsolution to the Cauchy problem~\eqref{Pi}; hence, the Comparison Principle ensures
\begin{equation*}
v^i(x)-k\leq w^i(t,x)\leq v^i(x)+k\qquad \forall (t,x)\in Q_\infty,
\end{equation*}
and, in particular: $\|w^i\|_\infty \leq 2k$.
Furthermore, by standard regularity theory for parabolic equations (see~\cite{CKS, W1,W2}) and by Proposition~\ref{Prp:erg}-($iii$), the function~$w^i$ fulfill hypothesis~\eqref{A100} with $\g=1$ and, consequently, also $u^i$ fulfill hypothesis~\eqref{A100} with $\g=1$.
\end{Proofc}

\begin{Proofc}{Proof of Theorem~\ref{thm:erg}: alternative version.}
We shall follow the arguments for the Comparison Principle (see \cite{CIL} and also \cite{JK2} for continuous dependence estimates).
For every positive $\eta$, define
\begin{equation}\label{psi_ell}
\psi(x,y):=w^1_\d(x)-w^2_\d(y)-\frac \eta 2 |x-y|^2
\end{equation}
where $w^i_\d$ ($i=1,2$) is the unique bounded (and periodic) solution to
\begin{equation}\label{Eidd}
\d w^i_\d+\min\limits_{\b\in B}\max\limits_{\a\in A}\left\{
-\tr\left(a_i(x,\a,\b)D^2w^i_\d\right)+f_i(x,\a,\b)\cdot Dw^i_\d +\ell_i(x,\a,\b)
\right\}=0\qquad \textrm{in }\R^n .
\end{equation}
(Here, taking advantage of the periodicity of $w^i_\d$, the penalization term is simpler than the one in the proof of Theorem~\ref{thm:dipcnt}.)
Owing to these properties of $w^i_\d$, we deduce that there exists a point~$(x_0,y_0)\in\R^n\times\R^n$ where the function~$\psi$ attains its global maximum.

Let us now claim that, for $\Co:=2K(1+\max\limits_{x,\a,\b,i}|\ell_i|)$ (where $K$ is the constant introduced in Proposition~\ref{Prp:erg}-($iv$)), there holds
\begin{equation}\label{stima:xy_ell}
\eta |x_0-y_0|\leq \Co.
\end{equation}
Actually, we observe that the inequality $\psi(x_0,x_0)+\psi(y_0,y_0)\leq 2\psi(x_0,y_0)$  gives
\begin{equation*}
\begin{array}{rcl}
\eta|x_0-y_0|^2&\leq& [w^1_\d(x_0)-w^1_\d(y_0)]+[w^2_\d(x_0)-w^2_\d(y_0)]\\
&\leq& 2K(1+\max\limits_{x,\a,\b,i}|\ell_i|)|x_0-y_0|
\end{array}
\end{equation*}
where the latter inequality is due to Proposition~\ref{Prp:erg}-($iv$); whence, estimate~\eqref{stima:xy_ell} easily follows.

By \cite[Theorem 3.2]{CIL}, for every $\nu>0$, there exist matrices $X,Y\in\Sym^n$ such that
\begin{eqnarray}\label{cil1_ell}
&\left(\eta(x_0-y_0),X\right)\in J^{2,+}w^1_\d(x_0),\qquad
\left(\eta(x_0-y_0),Y\right)\in J^{2,-}w^2_\d(y_0)&\\ \notag
&\left(\begin{array}{cc} X&0\\0&-Y\end{array}\right)\leq
\eta(1+2\nu\eta)\left(\begin{array}{cc} I&-I\\-I&I\end{array}\right).&
\end{eqnarray}
Moreover, by the same arguments as those in the proof of Theorem~\ref{thm:dipcnt}, for every $(\a,\b)\in A\times B$, from last inequality we deduce
\begin{equation*}
\tr\left(a_1(x_0,\a,\b)X\right)-\tr\left(a_2(y_0,\a,\b)Y\right)\leq
 \eta(1+2\nu\eta)\left\|\s_1(x_0,\a,\b)-\s_2(y_0,\a,\b)\right\|^2.
\end{equation*}
Since $w^1_\d$ (respectively, $w^2_\d$) is a subsolution (resp., a supersolution) to equation~\eqref{Eidd} with $i=1$ (resp., $i=2$), by relations~\eqref{cil1_ell}, we infer
\begin{eqnarray*}
\d w^1_\d(x_0)+\min\limits_{\b\in B}\max\limits_{\a\in A}\left\{
-\tr\left(a_1(x_0,\a,\b)X\right)+\eta f_1(x_0,\a,\b)\cdot (x_0-y_0) +\ell_1(x_0,\a,\b)
\right\}&\leq &0\\
\d w^2_\d(y_0)+\min\limits_{\b\in B}\max\limits_{\a\in A}\left\{
-\tr\left(a_2(y_0,\a,\b)Y\right)+\eta f_2(y_0,\a,\b)\cdot (x_0-y_0) +\ell_2(y_0,\a,\b) \right\}&\geq &0.
\end{eqnarray*}
Taking into account the last three inequalities, by the same calculations as before, we obtain
\begin{multline*}
\d \left( w^1_\d(x_0)- w^2_\d(y_0)\right)\leq
\eta(1+2\nu\eta)\max\limits_{\a,\b}\left\|\s_1(x_0,\a,\b)-\s_2(y_0,\a,\b)\right\|^2\\
+\eta|x_0-y_0|\max\limits_{\a,\b}\left|f_1(x_0,\a,\b)-f_2(y_0,\a,\b)\right|
+\max\limits_{\a,\b}\left|\ell_1(x_0,\a,\b)-\ell_2(y_0,\a,\b)\right|.
\end{multline*}
Letting $\nu\to 0^+$, by the regularity of the coefficients (see assumption ($A3$)) and by estimate~\eqref{stima:xy_ell}, 
we have
\begin{eqnarray}\notag
\d \left( w^1_\d(x_0)- w^2_\d(y_0)\right)&\leq&
2\eta\left(C^2_\s|x_0-y_0|^2+\max\limits_{x,\a,\b}\left\|\s_1-\s_2\right\|^2\right)
+\eta C_f|x_0-y_0|^2\\\notag &&\hskip 10mm
+\eta|x_0-y_0|\max\limits_{x,\a,\b}|f_1-f_2|
+\omega(|x_0-y_0|)+\max\limits_{x,\a,\b}|\ell_1-\ell_2|\\\notag &\leq&
\Co^2(2C_\s^2+C_f)\eta^{-1} +2\eta\max\limits_{x,\a,\b}\left\|\s_1-\s_2\right\|^2
+\Co\max\limits_{x,\a,\b}|f_1-f_2|
\\\label{finerg1}&&\hskip 10mm
 +\omega(\Co \eta^{-1})
+\max\limits_{x,\a,\b}|\ell_1-\ell_2|.
\end{eqnarray}
We separately consider the cases $\s_1=\s_2$ and $\s_1\neq \s_2$.
If $\s_1=\s_2$, as $\eta\to +\infty$, last inequality reads
\begin{equation*}
\d \left( w^1_\d(x_0)- w^2_\d(y_0)\right)\leq  
\Co\max\limits_{x,\a,\b}|f_1-f_2|+\max\limits_{x,\a,\b}|\ell_1-\ell_2|;
\end{equation*}
finally, as $\d\to 0^+$, we conclude
\begin{equation*}
U^1- U^2\leq  
\Co\max\limits_{x,\a,\b}|f_1-f_2|+\max\limits_{x,\a,\b}|\ell_1-\ell_2|.
\end{equation*}
If $\s_1\neq\s_2$, we choose $\eta=\Co/\max\limits_{x,\a,\b}\left\|\s_1-\s_2\right\|$; even though, in general, this is not the optimal choice for minimizing the right-hand side of~\eqref{finerg1}, the final estimate will behave with respect to $\Co$ in the desired manner for the purposes of section~\ref{pertsing}.
For $\tilde C:=2 C_\s^2+2+C_f$, we have:
\begin{multline*}
\d \left( w^1_\d(x_0)- w^2_\d(y_0)\right)\leq  
 \Co \left(\tilde C \max\limits_{x,\a,\b}\left\|\s_1-\s_2\right\|
+\max\limits_{x,\a,\b}|f_1-f_2|\right)
 +\omega(\max\limits_{x,\a,\b}\left\|\s_1-\s_2\right\|)\\
+\max\limits_{x,\a,\b}|\ell_1-\ell_2|;
\end{multline*}
finally, as $\d\to 0^+$, we conclude
\begin{equation*}
U^1- U^2\leq
 \Co \left(\tilde C\max\limits_{x,\a,\b}\left\|\s_1-\s_2\right\|
+\max\limits_{x,\a,\b}|f_1-f_2|\right)
 +\omega(\max\limits_{x,\a,\b}\left\|\s_1-\s_2\right\|)
+\max\limits_{x,\a,\b}|\ell_1-\ell_2|.
\end{equation*}
Hence, one side of the inequality of our statement is proved. Reversing the role of $w^1_\d$ and $w^2_\d$, one can easily obtain the other side; therefore, we shall omit its proof.
\end{Proofc}
\begin{Remark}\label{constant}
By the calculations of the proof above, a good choice is 
$\tilde M=2K(1+\max_{x,\a,\b,i}|\ell_i|)(2 C_\s^2+2+C_f)$, where $K$ is the constant introduced in Proposition~\ref{Prp:erg}-($iv$) while $C_\s$ and $C_f$ are the Lipschitz constants of~$\s$ and~$f$ respectively (see assumption~($A3$)).
\end{Remark}

%%%%%%%%%%%%%%%%%%%%%%%%%%%%%%%%
%%        Applicazione:       %%
%%   Perturbazioni singolari  %%
%%%%%%%%%%%%%%%%%%%%%%%%%%%%%%%%

\subsection{Singular perturbation problems}\label{pertsing}
We consider the following {\it singular perturbation} problems
\begin{equation}\label{SP}
\left\{
\begin{array}{ll}
\partial_t u^\e+{\cal H}\left(
x,y,D_xu^\e,\frac{D_yu^\e}\e,D^2_{xx}u^\e,\frac{D^2_{yy}u^\e}\e,\frac{D^2_{xy}u^\e}{\sqrt\e}\right)=0 &\quad\textrm{in }(0,T)\times\R^n\times\R^m\\
u^\e(0,x,y)=h(x) &\quad \textrm{on }\R^n\times\R^m
\end{array}\right.
\end{equation}
where $u^\e=u^\e(t,x,y)$ is a real function, $\e\in(0,1)$ and 
\begin{equation*}
{\cal H}(x,y,p,q,X,Y,Z):=\min\limits_{\b\in B}\max\limits_{\a\in A}
\left\{
-\tr(MX)-\tr(NY)-2\tr(EZ)+F\cdot q +G\cdot p +L
\right\}
\end{equation*}
with $\phi=\phi (x,y,\a,\b)$ for every $\phi=M,N,E,F,G,L$.
The aim of this section is to study the asymptotic behavior of $u^\e$ as $\e\to0^+$. For the wide literature on this matter, we refer the reader to the monographs by~Bensoussan~\cite{Ben}, Dontchev and~Zolezzi~\cite{DZ}, Kokotovi\'c, Khalil and O'Reilly~\cite{KKP}, Alvarez and Bardi~\cite{AB8} and references therein.
Let us only recall that these problems arise in zero-sum two-persons stochastic differential games \eqref{control}-\eqref{cost} where the state variable ``splits'' in the {\it slow} one $x$ and in the {\it fast} one $y$. For the control system
\begin{equation*}
\begin{array}{ll}
dx_s=G(x_s,y_s,\a_s,\b_s)+\sqrt 2 \Xi(x_s,y_s,\a_s,\b_s) dW_s,&\qquad x_0=x\\
dy_s=\e^{-1} f(x_s,y_s,\a_s,\b_s)+\sqrt{ 2\e^{-1}} \Sigma(x_s,y_s,\a_s,\b_s) dW_s,&\qquad y_0=y
\end{array}
\end{equation*}
and the cost functional
\begin{equation*}
P(t,x,y,\a,\b):=\mathbb E_{(x,y)}\left[ \int_0^tL(x_s,y_s,\a_s,\b_s)\, ds
+h(x_t,y_t)\right],
\end{equation*}
the lower value function $u^\e$ is a viscosity solution to problem~\eqref{SP} with $M={}^T\!\Xi\Xi$, $N={}^T\!\Sigma\Sigma$ and $E={}^T\!\Sigma\Xi$.

Throughout this section, we shall assume:
\begin{itemize}
\item[($S1$)] $A$ and $B$ are two compact metric spaces.
\item[($S2$)] $M={}^T\!\Xi\Xi$, $N={}^T\!\Sigma\Sigma$, $E={}^T\!\Sigma\Xi$. The functions $\Xi$, $\Sigma$, $F$, $G$, $L$ and~$h$ are bounded continuous functions in $\Rn\times\Rm\times A \times B$ with values respectively in~ $\M^{n,p}$, $\M^{m,p}$, $\Rm$, $\Rn$, $\R$ and $\R$, namely, for some $C>0$, there holds:
$\|\phi\|_\infty \leq C$ for $\phi=M,N,E,F,G,L$.

All these functions are $\mathbb Z^m$-periodic in $y$.
\item[($S3$)] The functions $\Xi$, $\Sigma$, $F$ and $G$ (respectively, $L$ and~$h$) are Lipschitz (resp., uniformly) continuous in $(x,y)$ uniformly in $(\a,\b)$
that is: there exists a positive constant $C_\phi$ and a modulus of continuity $\omega_\psi$ such that
\begin{eqnarray*}
\left|\phi(x_1,y_1,\a,\b)-\phi(x_2,y_2,\a,\b)\right|&\leq& C_\phi( |x_1-x_2|+|y_1-y_2|)\\
\left|\psi(x_1,y_1,\a,\b)-\psi(x_2,y_2,\a,\b)\right|&\leq& \omega_\psi(|x_1-x_2|+|y_1-y_2|)
\end{eqnarray*}
for every~$(x_i,y_i)\in\Rn\times\Rm$ ($i=1,2$) and~$(\a,\b)\in A\times B$, with $\phi=\Xi,\Sigma,F,G$ and $\psi=L,h$.
\item[($S4$)] There exists  $\nu>0$ such that, for every $(x,y,\a,\b)\in\Rn\times\Rm\times A\times B$, there holds
\begin{equation*}
M(x,y,\a,\b)\geq \nu I,\qquad N(x,y,\a,\b)\geq \nu I.
\end{equation*}
\end{itemize}

Let us recall from~\cite{AB2,AB8} the definition of the {\it effective} Hamiltonian~$\Ho$: for every $(\xo,\po,\Xo)\in\Rn\times\Rn\times\Sym^n$ fixed, the value~$-\Ho(\xo,\po,\Xo)$ is the ergodic constant for ${\cal H}(\xo,y,\po,q,\Xo,Y,0)$ with respect to the variable~$y$.
In other words, for $\d>0$, the problem
\begin{equation}\label{d-cp}
\d w_\d+{\cal H}\left(\xo,y,\po,D_y w,\Xo,D^2_{yy}w,0\right)=0
\qquad\textrm{in }\Rm,\qquad w_\d=w_\d(y)\quad\textrm{periodic}
\end{equation}
admits exactly one continuous solution and moreover, as $\d\to0^+$, $\d w_\d$ converges uniformly in $y$ to the value~$-\Ho(\xo,\po,\Xo)$. We refer the reader to Proposition~\ref{Prp:erg} for several properties of problem~\eqref{d-cp}.
In particular, let us observe (see also \cite[Theorem~4.1]{AB8}) that Proposition~\ref{Prp:erg}-($iv$) can be stated as follows: there exist $\kappa\in(0,1]$ and $K>0$ such that
\begin{equation}\label{pt-iv2}
\|w_\d -w_\d(0)\|_{C^{1,\kappa}}\leq K\left(1+|\po|+\|\Xo\|\right).
\end{equation}

%
%Proposizione
%

\begin{Proposition}\label{prp:pertsing}
The solution~$u^\e$ to problem \eqref{SP} converges locally uniformly in $[0,T)\times\Rn\times\Rm$ to the unique solution $u=u(t,x)$ to the effective problem
\begin{equation}\label{eff}
\left\{
\begin{array}{ll}
u_t+\Ho(x,D_xu,D^2_{xx}u)=0 &\qquad\textrm{in }(0,T)\times\Rn\\
u(0,x)=h(x)& \qquad \textrm{on }\R^n.
\end{array}\right.
\end{equation}
\end{Proposition}
\begin{Proofc}{Proof of Proposition~\ref{prp:pertsing}}
We shall argue using several results established by Alvarez and Bardi~\cite{AB2, AB8}.
Invoking~\cite[Theorem 2.9]{AB8} (see also~\cite[Corollary 2]{AB2}), it suffices to prove that the Comparison Principle applies to the effective problem~\eqref{eff}.
To this end, by virtue of the results by Ishii and Lions~\cite{il}, we need the following two properties: ($i$) $\Ho$ is uniformly elliptic, ($ii$) for some constant $\Ko$ and for some modulus of continuity $\bar\omega$, there holds
\begin{multline}\label{claim:CP}
\left|\Ho(x_1,p_1,X_1)-\Ho(x_2,p_2,X_2)\right|\leq
C\|X_1-X_2\| +C |p_1-p_2| +\bar\omega(|x_1-x_2|)\\
+\Ko |x_1-x_2|(1+|p_1|\vee|p_2|+\|X_1\|\vee\|X_2\|)
\end{multline}
for every $(x_i,p_i,X_i)\in \Rn\times\Rn\times\Sym^n$ ($i=1,2$).
We observe that the uniform ellipticity is well known so we shall omit its proof and we refer the reader to \cite[Theorem 4.4]{AB8} and \cite[Lemma 3.2]{Ev2} for the detailed proof.
In order to prove \eqref{claim:CP}, let us recall that, for $i=1,2$, the value $-\Ho(x_i,p_i,X_i)$ is the ergodic constant for the problem
\begin{multline*}
\min\limits_{\b\in B}\max\limits_{\a\in A}\left\{
-\tr(N(x_i,y,\a,\b)D^2_{yy}w^i)+D_y w^i\cdot F(x_i,y,\a,\b)
-\tr(M(x_i,y,\a,\b)X_i)\right.\\
\left.+p_i\cdot G(x_i,y,\a,\b)+L(x_i,y,\a,\b)\right\}
=-\Ho(x_i,p_i,X_i).
\end{multline*}
Applying Theorem~\ref{thm:erg} with the variable $x$ replaced by $y$ and
\begin{eqnarray*}
&\s_i(\cdot,\a,\b)=\Sigma(x_i,\cdot,\a,\b),\qquad
f_i(\cdot,\a,\b)=F(x_i,\cdot,\a,\b)&\\
& \ell_i(\cdot,\a,\b)=-\tr(M(x_i,\cdot,\a,\b)X_i)+p_i\cdot G(x_i,\cdot,\a,\b)+L(x_i,\cdot,\a,\b)&\\
&\omega(r)=\left[C_M (\|X_1\|\vee \|X_2\|)+C_G (|p_1|\vee |p_2|)\right]r+\omega_L(r)&
\end{eqnarray*}
for some constant $\tilde M$, we infer
\begin{align*}
&\left|\Ho(x_1,p_1,X_1)-\Ho(x_2,p_2,X_2)\right|\leq \\&\hskip 10mm
\tilde M
\left(\max\limits_{y,\a,\b}\left\|\Sigma(x_1,y,\a,\b)-\Sigma(x_2,y,\a,\b)\right\|
+\max\limits_{y,\a,\b}|F(x_1,y,\a,\b)-F(x_2,y,\a,\b)|\right)
\\&\hskip 10mm
+\left[C_M (\|X_1\|\vee \|X_2\|)+C_G (|p_1|\vee |p_2|)\right]
\max\limits_{y,\a,\b}\left\|\Sigma(x_1,y,\a,\b)-\Sigma(x_2,y,\a,\b)\right\|
 \\&\hskip 10mm
+\omega_L\left(\max\limits_{y,\a,\b}\left\|\Sigma(x_1,y,\a,\b)-\Sigma(x_2,y,\a,\b)\right\|\right)
 \\&\hskip 10mm
+\max\limits_{y,\a,\b}
\left|\tr\left[M(x_1,y,\a,\b)X_1-M(x_2,y,\a,\b)X_2\right]\right|
\\&\hskip 10mm
+\max\limits_{y,\a,\b}
\left|p_1\cdot G(x_1,y,\a,\b)-p_2\cdot G(x_2,y,\a,\b)\right|
+\max\limits_{y,\a,\b}
\left|L(x_1,y,\a,\b)-L(x_2,y,\a,\b)\right|.
\end{align*}
Taking into account the regularity of the coefficients (see assumption ($S3$)), we deduce
\begin{align*}
&\left|\Ho(x_1,p_1,X_1)-\Ho(x_2,p_2,X_2)\right|\leq
 \\&\hskip 10mm
C\|X_1-X_2\| +C|p_1-p_2|
+\omega_L(C_\Sigma |x_1-x_2|) 
+\omega_L\left(|x_1-x_2|\right)+\tilde M( C_\Sigma+C_F )|x_1-x_2|
\\&\hskip 10mm
+|x_1-x_2|\left[C_M (\|X_1\|\vee \|X_2\|)+C_G (|p_1|\vee |p_2|)\right]C_\Sigma 
+C_M |x_1-x_2|\left(\|X_1\|\wedge\|X_2\|\right)
\\&\hskip 10mm
+C_G |x_1-x_2|\left(|p_1|\wedge|p_2|\right).
\end{align*}
Since there holds $\max_{x,\a,\b,i}|\ell_i|\leq C\left(1+|p_1|\vee|p_2|+\|X_1\|\vee\|X_2\|\right)$, by Remark~\ref{constant}, we can choose
\begin{equation*}
\tilde M:=2K(C+1)
(2C^2_\Sigma+2+C_F)\left(1+|p_1|\vee|p_2|+\|X_1\|\vee\|X_2\|\right).
\end{equation*}
Hence, the previous inequality becomes
\begin{multline*}
\left|\Ho(x_1,p_1,X_1)-\Ho(x_2,p_2,X_2)\right|\leq
C\|X_1-X_2\| +C |p_1-p_2| +\omega_L(C_\Sigma|x_1-x_2|)+\omega_L(|x_1-x_2|)\\
+\Ko |x_1-x_2|(1+|p_1|\vee|p_2|+\|X_1\|\vee\|X_2\|)
\end{multline*}
for some constant $\Ko$ independent of $(x_i,p_i,X_i)$. Finally, choosing $\bar\omega(r):=\omega_L(C_\Sigma r)+\omega_L(r)$, our claim~\eqref{claim:CP} is completely proved.
\end{Proofc}

\end{document}